\documentclass[11pt]{amsart}

\usepackage{algorithm}
\usepackage{algpseudocode}
\usepackage[ansinew]{inputenc}
\usepackage[a4paper,dvips]{geometry}
\usepackage{latexsym}
\usepackage{amsfonts}
\usepackage{amsmath,amssymb,mathrsfs}
\usepackage{mathtools}
\usepackage{caption,subcaption}
\usepackage{bbm}
\usepackage[dvipsnames]{xcolor}
\usepackage{tikz}
\usetikzlibrary{decorations.pathreplacing,angles,quotes,calligraphy}

\usepackage{fixme}
\usepackage{graphicx}
\graphicspath{{images/}}
\usepackage{enumerate}
\usepackage{comment}

\usepackage[
backend=biber,
style=alphabetic,
sorting=ynt
]{biblatex}

\allowdisplaybreaks

\makeatletter
\newcommand{\addresseshere}{%
  \enddoc@text\let\enddoc@text\relax
}
\makeatother

\newtheorem{theorem}{Theorem}

\newtheorem{definition}[theorem]{Definition}

\newtheorem{conj}[theorem]{Conjecture}

\title[Improving the Halton sequence]{\large An improved Halton sequence for implementation in quasi-Monte Carlo methods}

\author{Nathan Kirk}
\address{Nathan Kirk \\ University of Waterloo, Ontario, Canada}
\email{n2kirk@uwaterloo.ca}

\author{Christiane Lemieux}
\address{Christiane Lemieux \\ University of Waterloo, Ontario, Canada}
\email{clemieux@uwaterloo.ca}

\date{\today}

\begin{document}

\maketitle

\begin{abstract}
Despite possessing the low-discrepancy property, the classical $d$-dim\-ensional Halton sequence is known to exhibit poorly distributed projections when $d$ becomes even moderately large. This, in turn, often implies bad performance when implemented in quasi-Monte Carlo (QMC) methods in comparison to, for example,  the Sobol' sequence. As an attempt to eradicate this issue, we propose an adapted Halton sequence built by integer \textit{and} irrational based van der Corput sequences and show empirically improved performance with respect to the accuracy of estimates in numerical integration and simulation. In addition, for the first time, a scrambling algorithm is proposed for irrational based digital sequences.
\end{abstract}

\section{INTRODUCTION}
\label{sec:intro}

\subsection{Context} Numerical integration is pivotal across various fields like computational mathematics, finance, physics and computer graphics where calculating integrals or expected values analytically often isn't feasible. For problems with one or a small number of variables, traditional integration techniques such as the trapezoidal rule are effective. However, contemporary challenges involve integrals dependent on a large number of variables implying a high dimensional integral, when it is necessary to employ more advanced methods, thus the development and improvement of numerical methods remains a vibrant field of study; \cite{MASK2024}, \cite{NUMERICALINT1}, \cite{NUMERICALINT2}, \cite{TRAUB2000}. In order to explain the contribution of this work, we introduce the problem formally. 

Given a finite set of sampling points $P_N = \{\mathbf{x}_1, \ldots, \mathbf{x}_N\}$ in the $d$-dimensional unit cube, $[0,1]^d$, one can approximate an integral by a standard so-called quadrature rule 
\begin{equation}\label{eq:approx}
I = \int_{[0,1]^d} f(\mathbf{x}) d\mathbf{x} \approx \frac{1}{N} \sum_{i=1}^N f(\mathbf{x}_i) =  \widehat{I}.
\end{equation}
Playing a crucial role throughout this text, an infinite sequence $P = \{\mathbf{x}_1, \mathbf{x}_2, \ldots \}$ contained in $[0,1]^d$ is said to be \textit{uniformly distributed modulo one} (abbreviated throughout as uniformly distributed) if and only if for every continuous real-valued function $f$ defined on $[0,1]^d$, we have $$ \lim_{N \rightarrow \infty} \frac{1}{N} \sum_{i=1}^N f(\mathbf{x}_i) = \int_{[0,1]^d} f(\mathbf{x}) d\mathbf{x}.$$ The quality of the approximation \eqref{eq:approx} can be bounded from above by the famous Koksma-Hlawka inequality \cite{HLAWKA1984}, \cite{KUIPNIED1974},
\begin{equation}\label{eq:KHinequal}
\left| I - \widehat{I} \right| \leq V(f)D(P_N).
\end{equation}
This result tells us that the absolute integration error can be bounded by the product of two factors. The first term $V(f)$ measures the variation of the integrand $f$ (in the sense of Hardy-Krause) and therefore depends on the problem at hand. The second, $D(P)$, is called the (star) discrepancy of $P$, and measures how uniformly the sampling points fill the space $[0,1]^d$. The more uniformly distributed the points are, the lower the discrepancy is and sequences with small enough discrepancy are called \textit{low-discrepancy}. One such standard procedure for approximating $I$, known as Monte Carlo (MC) integration, assumes that the sampling points are taken uniformly and randomly from $[0,1]^d$. As another possibility, one can attempt to ensure a `better' estimation by replacing the random sampling points with a low-discrepancy point set or sequence. This is then referred to as a quasi-Monte Carlo (QMC) method and such low-discrepancy point sets have been implemented in a wide range of disciplines, not solely in numerical integration; see \cite{GALANTIJUNG1997}, \cite{KELLER2022}, \cite{rusch2021} and \cite{TEYTAUD2020} for examples of QMC implementation. Fortunately, there exist many choices of low-discrepancy sets and sequences including $(t,m,s)$-nets, $(t,s)$-sequences, lattice rules and the subject of this work, the Halton sequence.

\subsection{The Halton Sequence}

As one of the first examples of a multidimensional low-discrepancy sequence introduced in \cite{HALTON1960}, the construction of the Halton sequence begins by considering one-dimensional van der Corput sequences in integer base $b \geq 2$.

\begin{definition}\label{def:intVDC}
Let $b \geq 2$ be an integer, then for $i \in \mathbb{N}_0$ we can write 
\[
i = \sum_{\ell=0}^{\infty} d_{i,\ell} b^{\ell}
\]
for \( d_{i, \ell} \in \mathbb{Z}_{b} \). Then the $i^{th}$ point of the \textit{van der Corput sequence in base $b$} is given by 
\[
g_{i,b} := \sum_{\ell=0}^\infty d_{i,\ell} b^{-(\ell+1)}.
\]
\end{definition}

Subsequently, the $d-$dimensional Halton sequence can be built by appending van der Corput sequences in base $b_j$ for $j \in \{1, \ldots, d\}$ in each coordinate where the bases are chosen to be any $d$ relatively prime natural numbers. In practice, the first $d$ primes are almost always used and in this text we will proceed with that assumption.

\begin{definition}
Let \( b_j \) be the \( j^{th} \) largest prime number for \( j \in \mathbb{N} \). For \( i \in \mathbb{N}_0 \), the $i^{th}$ point of the $d-$dimensional Halton sequence is defined as $$\left( g_{i,b_1}, g_{i,b_2}, 
\ldots, g_{i,b_d} \right).$$
\end{definition}

\begin{figure}[t]
  \centering
  \caption{Projection of the first $500$ points of the $26^{th}$ and $27^{th}$ of the classical Halton sequence (bases 101 and 103) and the interlaced Halton sequence (bases $19.462..$ and $20.05..$)}\label{fig:2627_classical}
  \begin{subfigure}[b]{0.43\textwidth}  % Reduced width of the subfigure
    \centering
    \includegraphics[width=4cm]{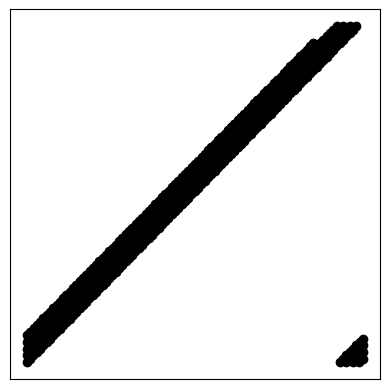}
    \caption*{}
    \label{fig:image1}
  \end{subfigure}
  \hspace{1mm}  % Adjust this space to bring images closer or further apart
  \begin{subfigure}[b]{0.43\textwidth}  % Reduced width of the subfigure
    \centering
    \includegraphics[width=4cm]{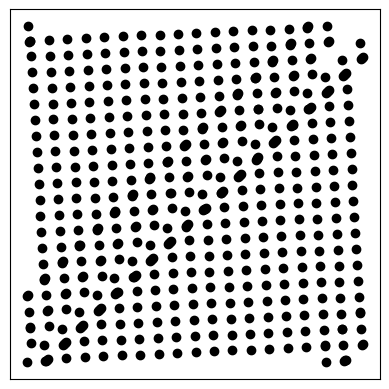}
    \caption*{}
    \label{fig:image2}
  \end{subfigure}

\end{figure}

The Halton sequence is known to be uniformly distributed; \cite{HALTON1960}, and even an example of a low-discrepancy sequence. However, it is also known that for moderate values of $N$ and larger $d$, the points tend to show unwanted structure; see for example, Figure \ref{fig:2627_classical} where the two dimensional projection of coordinates $26$ and $27$ (bases $101$ and $103$ respectively) of the Halton sequence   exhibit a structure which is somewhat co-linear. Larger dimensional coordinates, and in turn, larger bases, cause this unwanted distribution behaviour in the projections and in addition, is usually exacerbated by pairs of bases that are twin primes. There have been many efforts in recent decades to break this correlation between coordinates and improve the distribution properties of the Halton sequence. A popular approach, which has been thoroughly investigated in \cite{BRAATEN1979}, \cite{FAURE1992}, \cite{OKTEN2012}, \cite{OWEN2017}, \cite{VANDE2006}, \cite{HICKWANG2000}, is to induce some deterministic or random permutations on the digits of the integer base expansions which are appended to form the $d$ coordinates. 
In this work, the authors proceed in a different manner. Inspired by the recent work \cite{KIRKLEMIEUXWIART2023}, the authors introduce coordinates which are constructed from an irrational base to the classical Halton sequence in such a manner to delay the growth of the value of the van der Corput base, preserving good projected distribution properties for longer as the dimension increases.

\subsection{Previous Work}

Recently, \cite{KIRKLEMIEUXWIART2023} provided an extensive study of point sets and sequences constructed from an irrational base and introduced a complete framework analogous to that of classical $(t,m,s)-$nets and $(t,s)-$sequences in integer base $b$ to allow the study of equidistribution properties when working in an irrational base, $\gamma$, obtained as the largest root of $x^2-px-q$ for $1 \leq q \leq p$. In particular, as part of the same text, van der Corput sequences in irrational base $\gamma$ were defined in a similar number theoretic manner as in Definition \ref{def:intVDC} for the integer base case.

\begin{definition}[Definition 8.2, \cite{KIRKLEMIEUXWIART2023}]\label{def:irrVDC}
Let $\gamma$ be the largest root from $x^2-px-q$ for some choice of $1 \leq q \leq p$, and let $$n_i = \sum_{\ell=0}^\infty d_{i,\ell} (p+1)^{\ell}$$ be the digit expansion of the $i^{th}$ natural number which satisfies condition that any digit to the left of a $p$ in its base $(p+1)$ expansion is strictly less than $q$. Then the $i^{th}$ point of the van der Corput sequence (with parameters $(p,q)$) in base $\gamma$ is given by $$g_{n_i,\gamma} := \sum_{\ell=0}^\infty d_{i,\ell} \gamma^{-(\ell+1)}.$$
\end{definition}

Irrational based point sets and sequences were derived only in dimensions one and two in \cite{KIRKLEMIEUXWIART2023}, with the construction and study of higher dimensional constructions left as one of several open problems. Known to the authors, and referenced in their previous work, there had indeed been some small study of this multidimensional problem in terms of so-called $LS-$sequences. It is important to emphasise the similarities and differences between the two constructions. $LS-$sequences are a true generalization of the van der Corput sequences in integer base $b$, i.e., the integer sequences can be recovered by setting $L=b$ and $S=0$. The same cannot be said for our construction given precisely in Definition \ref{def:irrVDC}, however, it does have the advantage of being stated in an analogous manner to the classical integer sequences. Despite these differences, $LS-$sequences and van der Corput sequences in base $\gamma $ for $1 \leq q \leq p$ provide exactly the same construction when all parameters are strictly positive and when $L=p$ and $S=q$. 

$LS-$sequences were shown in \cite{CARBONE1} to be uniformly distributed in the unit interval if and only if $L \geq S$ and as a natural question, it was posed whether $LS-$sequences could be appended to create $d-$dimensional constructions similar to the generalization from van der Corput sequences to the classical Halton sequence in integer base. Multidimensional $LS-$sequences are investigated in \cite{AIS2012}, and it was shown (see Theorems 1 \& 2) that given parameters satisfying certain number theoretic conditions, these sequences not only fail to be uniformly distributed in $[0,1]^d$ but indeed, even fail to be dense in the space. Namely, re-written to match our notation, and assuming analogously to the Halton sequence that one-dimensional van der Corput sequences with parameters $(p,q)$ are simply appended coordinate-wise to form multidimensional sequences, the conditions under which the $d-$dimensional sequences fail to be uniformly distributed are stated below. The results below are written for two dimensions only, but of course, apply as a direct consequence to general $d \geq 2$ dimensions since a $d-$dimensional sequence is uniformly distributed only if all lower dimensional projections are also.

\vspace{2mm}
\noindent
\textbf{Condition A.} The two-dimensional sequence defined by the parameters $(p_1,q_1)$ and $(p_2,q_2)$ fails to be uniformly distributed in $[0,1]^2$ if there exists $k,m \in \mathbb{N}$ such that $\frac{\gamma_1^{k}}{\gamma_2^{m}} \in \mathbb{Q}$.

\vspace{2mm}
In words, this first condition implies that there cannot be any multiplicative relation between the irrational roots. For example, the two-dimensional sequence with parameters $(1,1)$ and $(11,1)$ fails to be uniformly distributed since $p_1=1,q_1=1 \implies \gamma_1 = \varphi$ and $p_2=11,q_2=1 \implies \gamma_2 = (11+5\sqrt{5})/2$, and we have $\gamma_1^5 = \gamma_2$. 

\vspace{2mm}
\noindent
\textbf{Condition B.} The two-dimensional sequence defined by the parameters $(p_1,q_1)$ and $(p_2,q_2)$ fails to be uniformly distributed in $[0,1]^2$ if $\text{gcd}(p_1,q_1,p_2,q_2) > 1$.

\vspace{2mm}
This second condition says that we not only have to worry about irrational-irrational interactions, but in fact also irrational-integer relations. For example, the two-dimensional sequence fails to be uniformly distributed when $p_1=2,q_1=0 \implies \gamma_1 = 2$ and $p_2=4,q_2=2 \implies \gamma_2 = 2+\sqrt{6}$ since $\text{gcd}(2,0,4,2) = 2 > 1$.

\vspace{2mm}
Stated explicitly in \cite{AIS2012}, ``up to date not a single example of parameters $(p_1,q_1,p_2,q_2)$ is known, for which either $q_1 \neq 0$ or $q_2 \neq 0$ and the corresponding two-dimensional sequence is uniformly distributed." In this paper, we show empirically that in fact, there are several good choices of parameters to generate $d-$dimensional constructions and we give one such algorithm to do so. It is important to note that the sequence that is presented has not been formally shown to be uniformly distributed, however, it does avoid the sufficient conditions A and B to not be uniformly distributed. In the next section, we state this as a formal conjecture to be proved. 
%Finally, these sequences massively improve upon the performance of the classical Halton sequence when implemented in a selection of quasi-Monte Carlo integration and simulation problems. 

\subsection{Paper Overview} Our paper is organized as follows. Section \ref{sec:interlaced_construction} introduces a new $d$-dimensional sequence, dubbed the interlaced Halton sequence, stating motivation for its construction and pseudo-code given in Algorithm \ref{algo:interlaced_construction}. A method for scrambling point sets and sequences constructed in irrational base (e.g., the van der Corput sequence in base $\gamma \in \mathbb{R}$) is given in Section \ref{sec:scrambling} with justification given by empirical evidence in the form of so-called $C_b(k; \cdot)$ values. Finally, Section \ref{sec:empirical_results} contains numerous numerical results from high-dimensional numerical integration and computational finance problems showing the positive impact the Halton sequence modification has in practice.

\section{THE INTERLACED HALTON SEQUENCE}

\subsection{The Construction}\label{sec:interlaced_construction}

\begin{algorithm}[h!]
    \caption{Generating Interlaced Halton Sequence}\label{algo:interlaced_construction}
    \begin{algorithmic}[1]
    \Statex \textbf{Input:} dimension $d$
    \Statex \textbf{Output:} list of bases $Bases$

    \State $Bases \leftarrow$ []

    \State $k \gets 0$ \Comment{$k$ is the index for the integer base in the classical Halton sequence}
    \State $p \gets 0$ \Comment{find valid $q$ for each $p$ and add it into suitable position}
    \While {len($Bases) < d$ }
    \State $p \gets p+1$
    \If{$p$ is prime}
    \State{$q \gets \left\lfloor \frac{p}{2} \right\rfloor$}
    \Else
    \State {$q = 1$}
    \EndIf
    
    \State $valid \gets False$
    
    \While{$q \leq p$}
    \State {$irbase \gets \text{larger root for } x^2 - px - q = 0$}

    \If{CheckPower($irbase$) \textbf{or} $\gcd(p, q) \neq 1$ \textbf{or} $\text{abs}(\text{round}(irbase, 1)-\text{round}(irbase, 0)) < 0.1$}
    \State \Comment{CheckPower is the function used to check Condition A}
    \State \Comment{i.e. if $irbase^k = base^m$ for some $base$ in the current sequence $Bases$}
    \State \Comment{We apply $k \in [2, 20]$ and $m \in [2, 20]$ in our algorithm}
    \State $q \gets q+1$
    \State \textbf{continue} \Comment{check conditions for $(p, q)$}
\EndIf

    \State $valid \gets True$
    \State \textbf{break}

    \EndWhile
    \State

    \If {not $valid$}
    \State continue
    \EndIf
    \Comment no valid $q$ for $p$, just ignore this $p$
    \State

    \If {$irbase >$ $k$-th integer base}
    \State $Bases.\text{append($k$-th integer base)}$ \Comment{add integer base from classical Halton sequence}
    \EndIf
    \State

    \If {len($Bases$) $== d$}
    \State \textbf{break}
    \EndIf
    \State

    \State $Bases.\text{append($irbase$)}$
   
    \EndWhile
    \State

    \State \textbf{Return} $Bases$    

\end{algorithmic}
\end{algorithm}
We now come to the primary contribution of this work as we introduce an algorithm to construct the $d-$dimensional \textit{interlaced Halton sequence}. When discussing the disadvantages of the classical Halton sequence earlier in the article, we noted that as the dimension increases, the quality of distribution in the projections becomes quite poor leading to inefficiency upon implementation in QMC methods. Thus, with a focus on mitigating this problem, the key idea of the construction presented in Algorithm \ref{algo:interlaced_construction} is to \textit{interlace} the Halton sequence by placing van der Corput sequences in an appropriate irrational base $\gamma$ between the existing integer base coordinates. 

%As a consequence, the rate at which the size of the base increases is reduced and therefore a greater level of uniformity can be maintained in the projections until a dimension much higher than that of the classical Halton sequence.

In words, assume $b_1$ and $b_2$ are consecutive prime numbers used as bases for subsequent coordinates of the classical Halton sequence. Then, for every natural number $b_1 \leq p < b_2$, Algorithm \ref{algo:interlaced_construction} builds the van der Corput sequences in base $\gamma$ with parameters $(p,1)$, unless $p$ is prime when we take $(p,\lfloor p/2 \rfloor)$, and these sequences are subsequently inserted in order of increasing $p$ between the base $b_1$ and $b_2$ van der Corput sequences. Of course, this is assuming the parameters chosen do not satisfy Conditions A or B. In the case that we are in breach of Conditions A or B, we incrementally increase the $q$ parameter until we find a suitable pair $(p,q)$, or if this is not possible, we ignore the current $p$ and move to the next natural number.

We make several comments on the above construction. Firstly, we have found that the probability of being in breach of Conditions A or B with a parameter pair $(p,q)$ is actually very small, particularly as $p$ increases. We also note that we have not attempted to provide `optimal' parameters for the interlaced Halton sequence. This search could cause computational issues, however, the most difficult obstacle to overcome may be that is it not obvious what criteria one should use for this search since the usual uniformity measures (such as the star or $L_2$-discrepancy) are known to be ineffective in very high dimensions; \cite{MORCAF1994}, \cite{WANGSLOAN2008}. Rather, Algorithm \ref{algo:interlaced_construction} can be implemented quickly to determine a $d-$dimensional sequence with more uniform projections that the classical Halton sequence. Finally, it should be made clear to the reader that several modifications of the above construction could be made. For example, given a pair of consecutive primes, $b_1$ and $b_2$, a more simple approach would be to interlace a single van der Corput sequence in base $\gamma$ where $\gamma \approx (b_1+b_2)/2$ while still avoiding conflict with Conditions A and B. After extensive numerical experimentation, we have found most success in numerical experimentation with the sequence as defined above and all results presented in Section \ref{sec:empirical_results} were completed with the construction in Algorithm \ref{algo:interlaced_construction}.

\begin{conj}\label{conj:UD}
    The $d-$dimensional interlaced Halton sequence in Algorithm \ref{algo:interlaced_construction} is uniformly distributed.
\end{conj}

When proving the uniformity of the classical integer Halton sequence, one is required to apply the famous Chinese remainder theorem. Hence, a formal proof of the above conjecture seems difficult in nature due to the introduction of irrationals to the Halton construction.

\subsection{Scrambling in Irrational Base}\label{sec:scrambling}

Despite the renewed interest in constructing digital nets and sequences from an irrational base, a major downside, as discussed in \cite{KIRKLEMIEUXWIART2023}, involves the difficulty of randomizing such sequences. When working with an integer base, randomization of a low-discrepancy point set allows researchers to easily derive error estimates when implementing QMC point sets in approximation methods which is then usually referred to in the literature as a randomized quasi-Monte Carlo (RQMC) method; see \cite{PIERRE2018} for further reading. There exist several popular methods for randomization with one of the most favored being so-called \textit{nested uniform scrambling} (or, Owen scrambling) introduced in \cite{OWEN1995}. Scrambling can enhance the uniformity of an already existing low-discrepancy point set, however, the attribute most useful for application purposes is the ability to subdue or completely remove any strong correlation between coordinates of a multidimensional point set. In this section, we provide a functional scrambling solution for point sets and sequences constructed in irrational base $\gamma$ and apply this algorithm to randomize the interlaced Halton sequence. This partially and empirically solves an open question from \cite{KIRKLEMIEUXWIART2023}. To motivate our approach, we first review key results from \cite{DONGLEMIEUXWIART2021} and \cite{DONGLEMIEUX2022} that allow us to characterize the potential for our construction to benefit from scrambling in the base we propose to use for scrambling. 

The quantity first introduced by \cite{DONGLEMIEUXWIART2021} and adapted to the case of the Halton sequence in \cite{DONGLEMIEUX2022} can be used to measure the above mentioned propensity of an arbitrary point set $P_N$ contained in $[0,1]^d$ to become very uniform after scrambling is called $C_b(\mathbf{k};P_N)$ values. In the case of a Halton sequence, the base $b$ is instead a vector $\mathbf{b}$ containing the integers used for scrambling each of the $d$ coordinates. Our construction inherits this feature of the Halton sequence, so henceforth we assume we are working with such vectors $\mathbf{b}$. To define this idea formally, we first have to explain the concept of an elementary interval in base $b$. For $b \geq 2$, an \emph{elementary $k$-interval in base $b$} is a subset of $[0,1)$ of the form
\[
\left[\frac{a}{b^k},\frac{a+1}{b^k}\right)
\]
where $0\leq a<b^k$ and $k \in \mathbb{N}_0$. In $d-$dimensions, an \emph{elementary $(k_1,\dots,k_s)-$interval in (multi-base) $\mathbf{b} = (b_1,\ldots,b_k)$} is the product of elementary $k_j-$intervals in base $b_j$, i.e. it is a set of the form
\[
\prod_{j=1}^d\left[\frac{a_j}{b_j^{k_j}},\frac{a_j+1}{b_j^{k_j}}\right),
\]
where $0\leq a_j<b_j^{k_j}$ for each $1\leq j \leq d$. The $C_{\mathbf{b}}(\mathbf{k}; P_N)$ values are intrinsically linked to how many distinct pairs of points from $P_N$ lie in the same elementary $\mathbf{k}$-interval in base $\mathbf{b}$.

\begin{definition}[Definition 4.1, \cite{DONGLEMIEUX2022}]
    Let $P_N$ be a set of $N$ points in $[0,1]^d$ and $b\geq 2$ be an integer. Let $C_{\mathbf{b}}(\mathbf{k};P_N)$ be defined as $$C_{\mathbf{b}}(\mathbf{k};P_N) := \frac{\prod_{j=1}^d b_j^{k_j}  M_{\mathbf{b}}(\mathbf{k};P_N)}{N(N-1)}$$ where $M_{\mathbf{b}}(\mathbf{k};P_N)$ denotes the number of ordered pairs of distinct points from $P_N$ that lie in the same elementary $\mathbf{k}-$intervals. 
\end{definition}

Being able to compute the $C_{\mathbf{b}}(\mathbf{k};P_N)$ values for a given point set informs the user about the distribution of $P_N$ in a very informative manner and importantly, this notion does not require the point set $P_N$ to hold any particular structure, nor in fact does it place any restriction on the number of points in $P_N$.

\begin{definition}[Definition 4.4, \cite{DONGLEMIEUX2022}]
Let $P_N$ be a point set of size $N$ in $[0, 1]^d$ and $b_j \geq 2,j=1,\ldots,d$ be $d$ bases. Let $\mathbf{k} = (k_1, \ldots, k_d) \in \mathbb{N}_0^d.$ Then we say $P_N$ is $\mathbf{k}-$quasi-equidistributed in base $\mathbf{b}$ if $C_{\mathbf{b}}(\mathbf{k}; P_N) \leq 1$. If $C_{\mathbf{b}}(\mathbf{k}; P_N) \leq 1$ for all $\mathbf{k} \in \mathbb{N}_0^d$ then we say $P_N$ is completely quasi-equidistributed (c.q.e) in base $\mathbf{b}$.    
\end{definition}

A key property, widely discussed in \cite{DONGLEMIEUXWIART2021}, is the propensity of a point set $P_N$ to be highly uniform after scrambling in base $b$ if the $C_b(\mathbf{k};P_N)$ values are ``small", i.e. if $P_N$ is completely quasi-equidistributed in base $b$. The same idea applies to a point set---such as the Halton sequence---for which a scrambling in base $b_j$ is applied to the $j$-th coordinate.

For our interlaced Halton sequence, we propose to scramble the coordinates associated with an irrational base, characterized by the pair $(p,q)$, by choosing $b=p+1$ as the scrambling base. Our proposed approach is based on the fact that
after much numerical investigation, we came up with the following conjecture.

\begin{conj}\label{conj:CQE}
    Let $\gamma$ be the largest root of the polynomial $x^2-px-q$ for $1 \leq q \leq p$. The van der Corput sequence in base $\gamma$ is c.q.e. in base $(p+1)$.
\end{conj}

We give one example of the golden ratio van der Corput sequence (with parameters $p=q=1$) below and remark that similar results were observed independent of the choice of $p$ and $q$ while scrambling in base $p+1$. Assuming Conjecture \ref{conj:CQE} holds, the van der Corput sequence in base $\gamma$ obtained from parameters $(p,q)$ can be scrambled in base $(p+1)$ to facilitate the desired randomization. 

\begin{table}[htb]
\centering
\caption{The set of the first $500$ points of the van der Corput sequence in base golden ratio is completely quasi-equidistributed base $2$.}
\begin{tabular}{|c||*{9}{c|}}
\hline
$k$ & 1 & 2 & 3 & 4 & 5 & 6 & 7 & 8 & 9 \\ \hline
$C_2(k; \cdot)$ & $0.9980$ & $0.9940$ & $0.9861$ & $0.9706$ & $0.9393$ & $0.8854$ & $0.7726$ & $0.6074$ & $0.2257$ \\ \hline
\end{tabular}
\label{your-table-label}
\end{table}

\section{Empirical Results}\label{sec:empirical_results}

\subsection{Numerical Integration} To test the distribution of sequences and point sets in the unit cube, one usually computes the star or $L_2$-discrepancy. In higher dimensions, however, this can be extremely computationally intensive and in addition, sometimes does not tell the full story; see \cite{WANGSLOAN2008}. As an alternative, suggested by the Koksma-Hlawka inequality \eqref{eq:KHinequal}, one can implement uniformly distributed point sets in quasi-Monte Carlo numerical integration of various standard test functions to compare the error of the approximation. The intuition here is that the more uniform the point set is, the faster the error convergence will be. In this section, we test the interlaced Halton sequence against the classical Halton in several test functions. The first of which, first proposed in \cite{SOBOL1996} and commonly used in similar studies such as \cite{OWEN2003} and \cite{LEMIEUXFAURE2009} is $$f_1(\mathbf{x}) = \prod_{j=1}^d \frac{\left| 4x_j -2 \right| + a_j}{1 + a_j}$$ for $\mathbf{x} = (x_1,\ldots, x_d) \in [0,1]^d$ and coefficient vector $\mathbf{a} = (a_1,\ldots,a_d) \in \mathbb{R}^d$. We give two variants of the coefficient vector $\mathbf{a}$ which dictates the \textit{effective dimension} of the problem; we refer to \cite{CAFL1997} for the full details on the effective dimension and note for now that it can be informally described as the number of dimensions which explain most of the problem's variance.
%At a high level, the effective dimension of a problem is the number of variables which explains 99\% of a problem's variance. It is typically less than the nominal one, and a smaller effective dimension often means the function can be more easily approximated through numerical methods; we refer to \cite{CAFL1997} for the details.
%We will test case (1) $a_j = j^2$ and case (2) $a_j = j$. In these two cases, the \textit{effective dimension} is small and large respectively, i.e., the integral with case (2) is difficult to compute as variables remain of equal importance with increasing dimension, contrary to case (1) where the variables become less and less important for the integral value.

As a second test function, previously used in \cite{SOBOL2003}, we take $$f_2(\mathbf{x}) = \prod_{j=1}^d 1 + c (x_j - 0.5)$$ where $c\in \mathbb{R}$. For this problem, the effective dimension is determined by the pair of parameters $(c,d)$ and again, we explore different combinations of parameters. 

%We will consider again two cases, (1) the pair $(0.1,120)$ and (2) the pair $(1,150)$ where is has been calculated that the effective dimension is $4$ and $20$ respectively.

\begin{comment}
As a final test function, we consider the problem of integrating the Keister function with respect to a $d-$dimensional Gaussian measure. That is,
\begin{equation*}
f_3(x) = \pi^{d/2} \cos(\lVert \mathbf{x} \rVert), \quad \mathbf{x} \in \mathbb{R}^d, \quad \mathbf{X} \sim N(0_d, I_d/2),
\end{equation*}

\begin{equation*}
\mu = \mathbb{E}[f_3(\mathbf{X})] := \int_{\mathbb{R}^d} f_3(x) \pi^{-d/2} \exp(-\lVert x \rVert^2) \, dx
\end{equation*}

\begin{equation}\label{eq:integral3}
= \int_{[0,1]^d} \pi^{d/2} \cos\left( \sqrt{\frac{1}{2} \sum_{j=1}^{d} \Phi^{-1}(x_j)} \right) \, dx
\end{equation}

The integrals of $f_1$ and $f_2$ over the $d$-dimensional unit cube both evaluate to $1$.
\end{comment}

\subsection{Simulation in Computational Finance}

Comparing quasi-Monte Carlo point sets via toy test function integrands can be viewed as subjective. Therefore, researchers often use integrals appearing in real-world situations to assess the performance comparison. In this text, we utilize a common application in computational finance and attempt to estimate the value at time $0$ of an Asian call option on an underlying asset that follows a log-normal distribution. The value at time $0$ of this option is given by the expectation, $$C_0 = \mathbb{E} \left[ e^{-rT} \left( \frac{1}{d} \sum_{j=1}^{d} S(u_j) - K \right)^+ \right]$$ where we let \( T \) be the expiration time of the contract, \( K \) the strike price, $S(u)$ 
%= S(0)e^{\left(r-\frac{\sigma^2}{2}\right)u+\sigma\sqrt{u}Z}\) where $Z \sim N(0,1)$ 
is the price of the underlying asset at time \( u \), and \( 0 < u_1 < \ldots < u_d = T \) are \( d \) times at which the asset price is observed, and the expectation is taken under the risk-neutral probability measure. Finally, \( r \) is the risk-free rate, the notation \( x^+ \) means \( \max(0, x) \), $\Phi^{-1}$ is the inverse CDF of the standard normal distribution and $\Delta_l = u_l - u_{l-1}$. Assuming the stock price follows a geometric Brownian motion with volatility $\sigma$, it can be shown that this expectation can be written as follows:
\begin{equation}
C_0 = e^{-rT} \int_{[0,1]^d} \left( \frac{1}{d} \sum_{j=1}^{d} S(0)e^{(r-\frac{\sigma^2}{2})u_j+\sigma  \sum_{l=1}^{j} \sqrt{\Delta_l} \Phi^{-1}(x_l)} - K \right)^+ \, dx_1 \ldots dx_d.
\end{equation}

We will proceed with the parameters $S(0) = 50, T = 1 \text{ year }, r = 0.05, \sigma = 0.3, K=45$ and consider $d=50$ and $d=100$. The true value of $C_0$ for this given set of parameters is calculated in advance of numerical approximation comparisons via QMC integration with $2 \times 10^6$ Sobol' points.

\subsection{Results} 
\begin{figure}[htb]
\centering
% Row 1
\caption{Integration error of function $f_1$. (\textbf{Top}) $a_j = j$ (\textbf{Bottom}) $a_j = j^2$.}\label{fig:f1_results}
\begin{subfigure}{0.33\textwidth}
    \includegraphics[width=\linewidth]{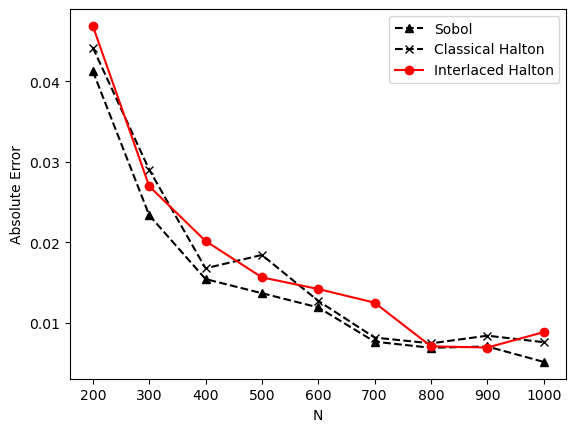}
    \caption*{$d=25$} % empty caption
\end{subfigure}\hfill
\begin{subfigure}{0.33\textwidth}
    \includegraphics[width=\linewidth]{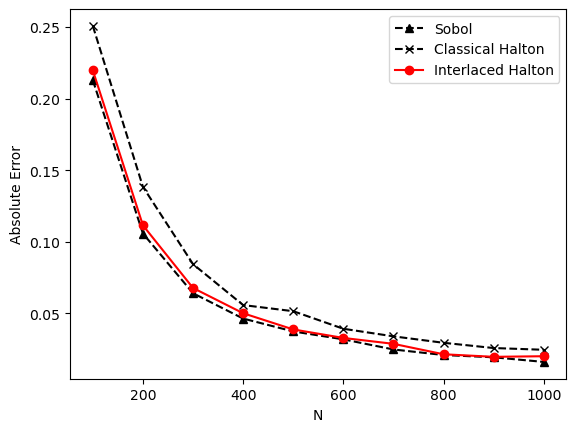}
    \caption*{$d=50$}
\end{subfigure}\hfill
\begin{subfigure}{0.33\textwidth}
    \includegraphics[width=\linewidth]{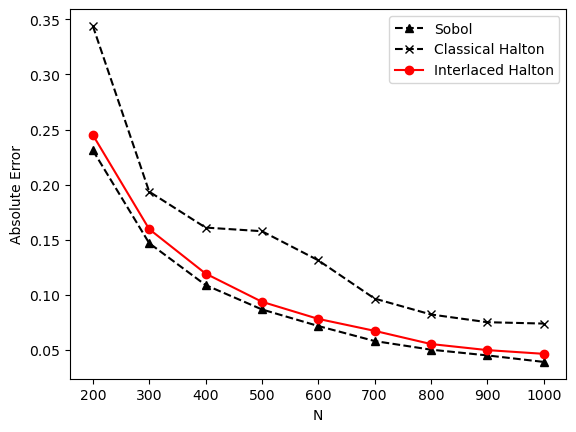}
    \caption*{$d=100$}
\end{subfigure}

% Row 2
\begin{subfigure}{0.33\textwidth}
    \includegraphics[width=\linewidth]{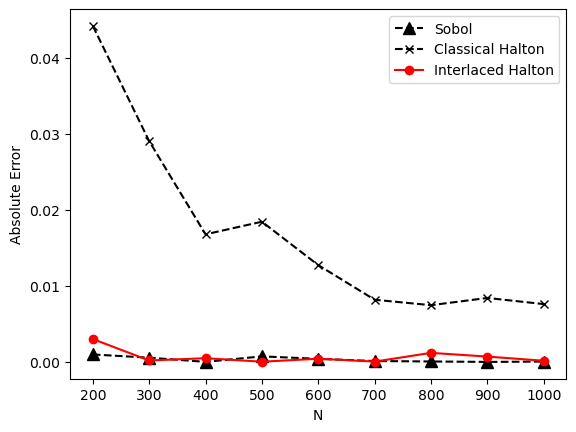}
    \caption*{}
\end{subfigure}\hfill
\begin{subfigure}{0.33\textwidth}
    \includegraphics[width=\linewidth]{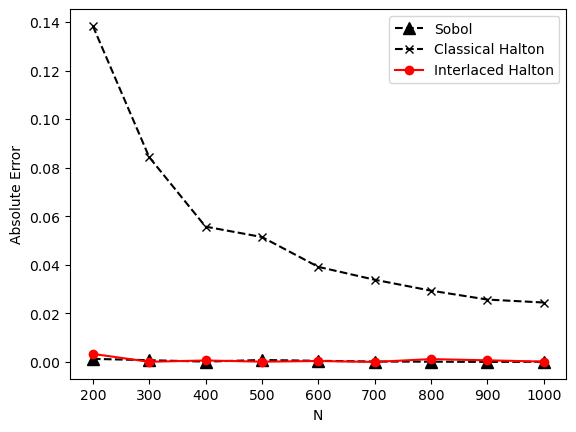}
    \caption*{}
\end{subfigure}\hfill
\begin{subfigure}{0.33\textwidth}
    \includegraphics[width=\linewidth]{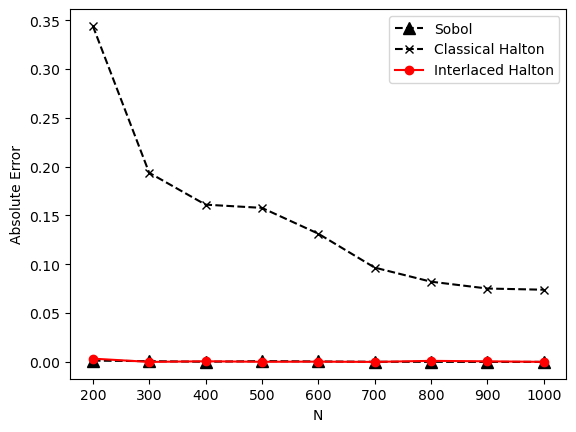}
    \caption*{}
\end{subfigure}
\end{figure}

\begin{figure}[htb]
\centering
% Row 1
\caption{Integration error of function $f_2$. (\textbf{Top}) $c=1$ (\textbf{Bottom}) $c=0.1$.}\label{fig:f2_results}
\begin{subfigure}{0.33\textwidth}
    \includegraphics[width=\linewidth]{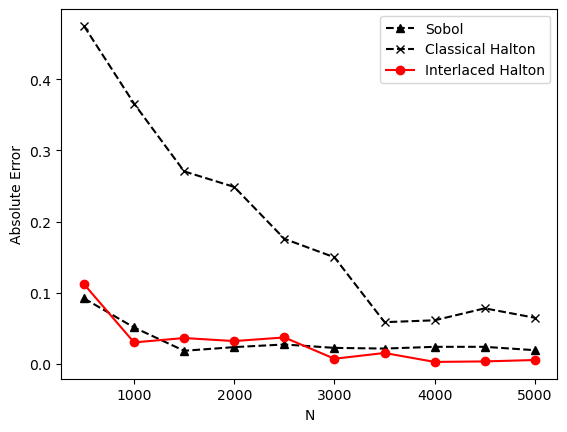}
    \caption*{$d=25$}
\end{subfigure}\hfill
\begin{subfigure}{0.33\textwidth}
    \includegraphics[width=\linewidth]{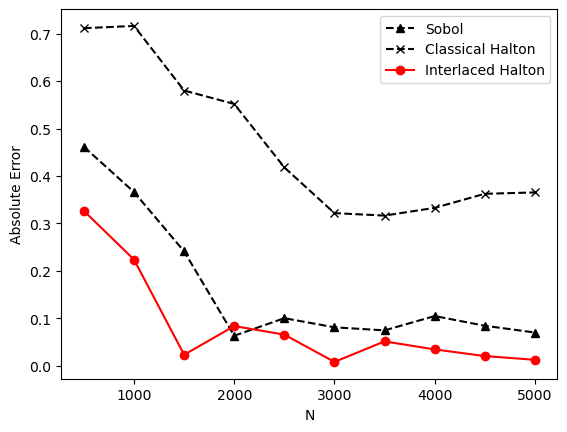}
    \caption*{$d=50$}
\end{subfigure}\hfill
\begin{subfigure}{0.33\textwidth}
    \includegraphics[width=\linewidth]{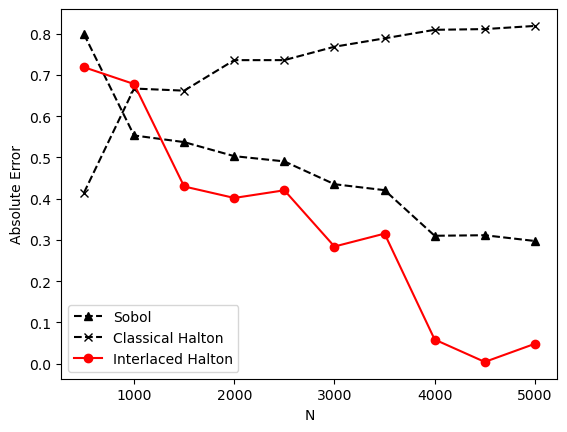}
    \caption*{$d=100$}
\end{subfigure}

\begin{subfigure}{0.33\textwidth}
    \includegraphics[width=\linewidth]{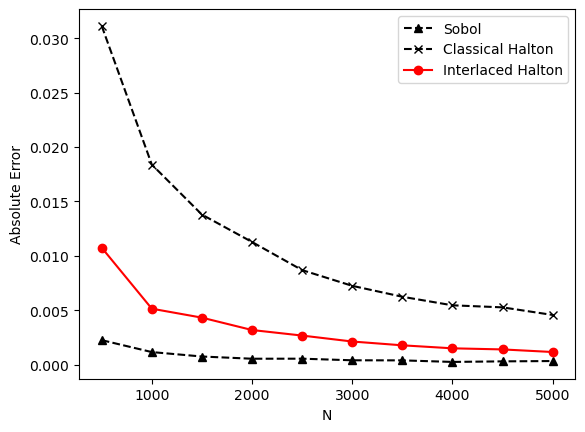}
    \caption*{} % empty caption
\end{subfigure}\hfill
\begin{subfigure}{0.33\textwidth}
    \includegraphics[width=\linewidth]{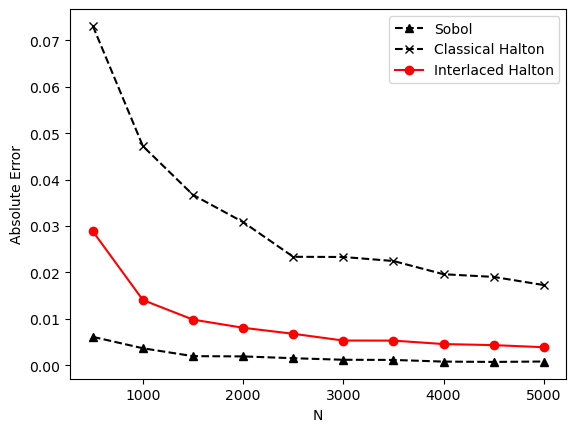}
    \caption*{}
\end{subfigure}\hfill
\begin{subfigure}{0.33\textwidth}
    \includegraphics[width=\linewidth]{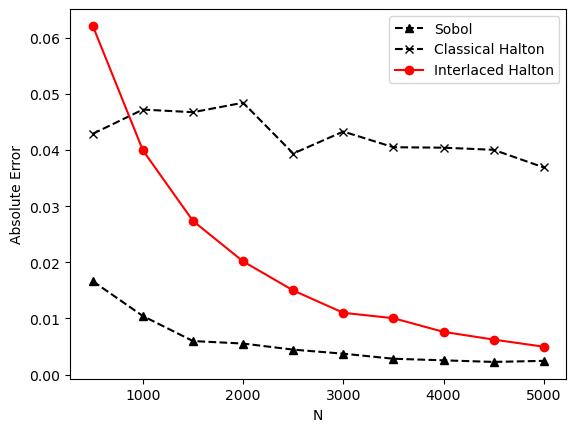}
    \caption*{}
\end{subfigure}

\end{figure}

For each of the functions $f_1$ and $f_2$, we have considered the integral approximations in dimensions $d=25, 50$ and $100$, and plotted the absolute (deterministic) error for each function for ranges of $N$ which best illustrate the error comparison. To further diversify our numerical results, for each integral, we have selected two choices of parameters to reflect differences in the problem's effective dimension. The results from quadrature approximations can be found in Figures \ref{fig:f1_results} and \ref{fig:f2_results}. We point out that each graph has its own error scale on the vertical axis. Further, for both $f_1$ and $f_2$ the parameters chosen and displayed in the top lines of each figure are known to be of higher effective dimension and are therefore more difficult to approximate accurately than the bottom line. For the computational finance problem, we considered both a deterministic error and the estimated variance (from 20 replications) derived from implementing randomized versions of each sequence. The Sobol' and classical Halton sequence were scrambled via the \texttt{qmcpy} while the interlaced Halton was scrambled by computer code implementing nested integer based scrambling.

\begin{figure}[htb]
\centering
\caption{Absolute error (top) and estimated variance (bottom) for Asian option problem. (\textbf{Left}) $d=50$ (\textbf{Right}) $d=100$.}\label{fig:asian_results}
\begin{subfigure}{0.4\textwidth}
    \includegraphics[width=\linewidth]{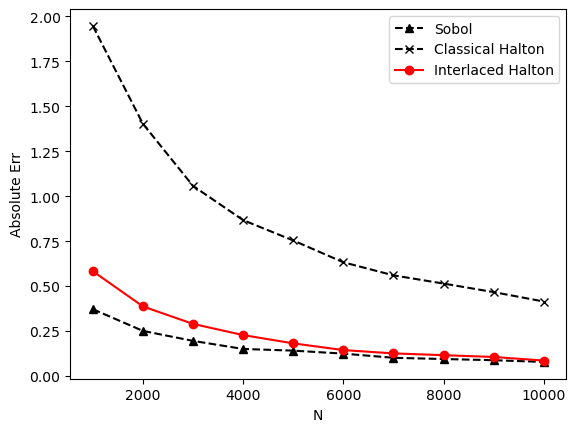}
    \caption*{}
\end{subfigure}
\begin{subfigure}{0.4\textwidth}
    \includegraphics[width=\linewidth]{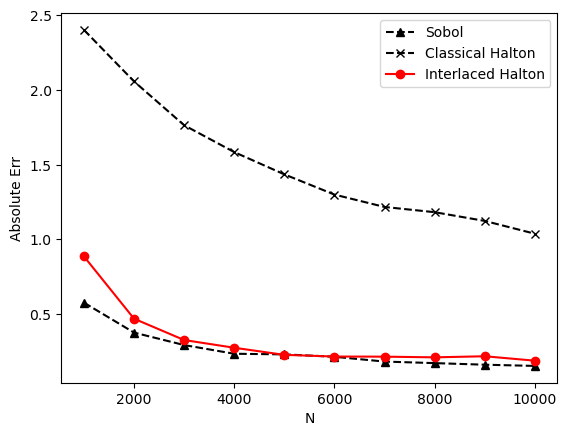}
    \caption*{}
\end{subfigure}

\begin{subfigure}{0.4\textwidth}
    \includegraphics[width=\linewidth]{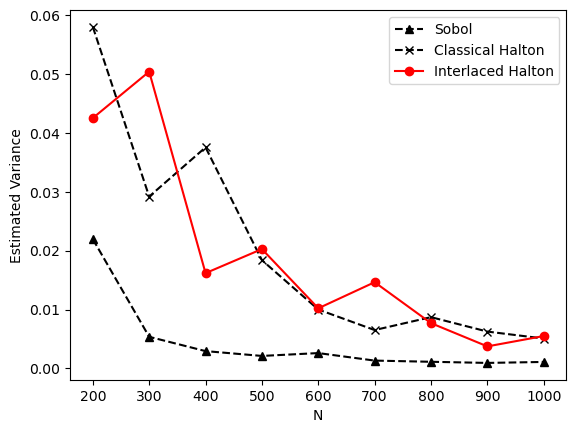}
    \caption*{} % empty caption
\end{subfigure}
\begin{subfigure}{0.4\textwidth}
    \includegraphics[width=\linewidth]{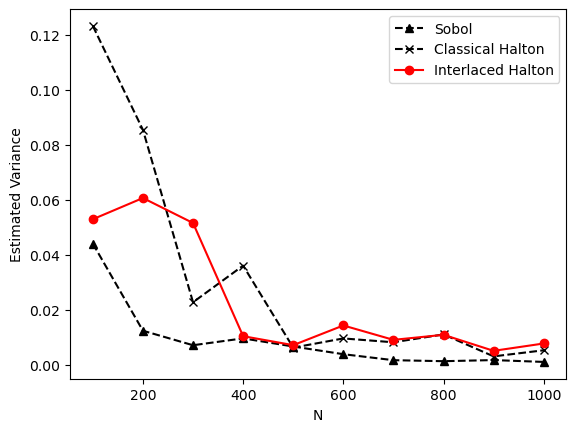}
    \caption*{}
\end{subfigure}
\end{figure}

The authors wish to highlight several important features which are the primary inspiration of this work. Firstly, across all experiments, the interlaced Halton sequence provides a smaller error in approximations and simulations than the classical Halton sequence. Specifically, it seems reasonable to conclude that, as the dimension increases, the interlaced Halton sequence appears to perform better than the classical Halton sequence by a greater margin and even, oftentimes, provides approximations that are competitive with that of the Sobol' sequence. The authors claim that the observed improvements are attributed to the enhancements in the lower-dimensional projections of the interlaced Halton sequence. As last comment, the small variance results contained in Figure \ref{fig:asian_results} empirically validate that scrambling the interlaced Halton with integer bases is a functional and effective approach to randomization of point sets and sequences constructed with an irrational base.

\subsection*{Future Work} We pose the following directions for future work. In the first instance, it would be of primary interest to formally prove Conjectures \ref{conj:UD} and \ref{conj:CQE} to verify the uniform distribution property of the interlaced Halton sequence and provide a solid foundation to encourage future study. It would also be interesting to determine which parameters of $p$ and $q$ generate the most improved performance over the classical Halton sequence in numerical integration or simulation problems. Finally, while we have shown that scrambling irrational based sequences in an integer base is a practical and functional option for randomization, the possibility of scrambling in irrational base $\gamma$ while retaining the desired properties still remains unknown.

\section*{ACKNOWLEDGMENTS}
The authors would like to thank Melody Mao for her assistance with numerical computations during her time as an undergraduate student at the University of Waterloo. We also acknowledge the use of computer code by Gracia Dong (University of Toronto) for calculations throughout. This work was supported by the Natural Science and Engineering Research Council of Canada (NSERC) via grant 238959.

% Reducing font size (to 9pt) for References & Author Biagraphies
\footnotesize

% AUTHOR: Include your bib file here
\printbibliography

\end{document}